\documentclass[leqno,12pt]{article}
\usepackage{amssymb,amsfonts,siunitx}
\usepackage{amsmath,latexsym}

\usepackage{amstext,epic,eepic,epsf,pslatex}
\usepackage{graphicx}

\newcommand{\R}{{\mathbb R}}

\newcommand{\Tr}{\hbox{Tr\,}}

\def\XXint#1#2#3{{\setbox0=\hbox{$#1{#2#3}{\int}$ }
\vcenter{\hbox{$#2#3$ }}\kern-.6\wd0}}

\newtheorem{defin}{Definition}[section] 
 
\newtheorem{thm}{Theorem}[section]

\begin{document}

\title{Symmetric Divergence-free tensors in the Calculus of Variations}

\author{Denis Serre \\ \'Ecole Normale Sup\'erieure de Lyon\thanks{U.M.P.A., UMR CNRS--ENSL \# 5669. 46 all\'ee d'Italie, 69364 Lyon cedex 07. France. {\tt denis.serre@ens-lyon.fr}}}

\date{}

\maketitle

\begin{abstract}
Divergence-free  symmetric tensors seem ubiquitous in Mathematical Physics. We show that this structure occurs in models that are described by the so-called ``second'' variational principle, where the argument of the Lagrangian is a closed differential form. Divergence-free tensors are nothing but the second form of the Euler--Lagrange equations. The symmetry is associated with the invariance of the Lagrangian density upon the action of some orthogonal group.
\end{abstract}


\paragraph{Notations.} The $n\times n$ identity matrix is $I_n$. Given two vectors $X,Y$, one denotes $X\otimes Y=XY^T$ the rank-one matrix with entries $x_iy_j$. With a capital letter D, the operator Div applies row-wise to matrix-valued fields (tensors)~:
$$({\rm Div}\,S)_i=\sum_j\partial_js_{ij}.$$
The Divergence of a tensor is thus a vector field.
A tensor $S$ is {\em Divergence-free} if ${\rm Div}\,S\equiv0$ in the sense of distributions.

\section{Ordinary Calculus of Variations}\label{s:clas}

In its classical form, Calculus of Variations (CV) deals with a functional
$${\cal F}[u]:=\int_D L(u,\nabla u)\,dy$$
where $D$ is a $d$-dimensional open domain, the argument $u:D\rightarrow\R$ is a function of some Sobolev class (ensuring that ${\cal F}[u]$ is well-defined) and $L:\R^{1+d}\rightarrow\R$ is a smooth given function. For instance, $L(u,p)=\ell(p)=O(|p|^m)$ and $u$ is taken in $W^{1,m}(D)$.

A Variational Principle (VP) declares that $u$ is a critical point if $\delta{\cal F}[u]=0$, where $\delta$ is differentiation along admissible paths passing through $u$. Paths $\epsilon\mapsto v_\epsilon$ must be such that $v_0=u$ and  $\epsilon\mapsto{\cal F}[v_\epsilon]$ is differentiable. We ask in addition that $v_\epsilon\equiv u$ away from a compact subset of $D$. Two candidates emerge naturally, the first being $v_\epsilon=u+\epsilon h$ for some $h\in{\cal D}(D)$. The second is the little more involved $v_\epsilon=u\circ\phi_\epsilon$ where $(\phi_\epsilon)_{\epsilon\in\R}$ is the flow of some vector field $\xi\in{\cal D}(D)^d$~:
$$\frac{d \phi_\epsilon}{d\epsilon}=\xi(\phi_\epsilon),\qquad\phi_0={\rm Id}_D.$$
The corresponding conditions
$$\left.\frac{d{\cal F}[v_\epsilon]}{d\epsilon}\right|_{\epsilon=0}=0,\qquad\forall\, h\hbox{ or }\forall\, \xi$$
are called the {\em first} and {\em second formulations} of the VP, respectively. They are not always equivalent to each other. The first one is well-defined for reasonable Lagrangian densities, for instance it suffices that $\frac{\partial L}{\partial p}(u,\nabla u)$ and  $\frac{\partial L}{\partial u}(u,\nabla u)$  be uniformly locally 
integrable whenever $u\in W^{1,m}(D)$. It yields the well-known {\em first form} of the Euler--Lagrange equations
\begin{equation}
\label{eq:ffEL}
{\rm div}\frac{\partial L}{\partial p}=\frac{\partial L}{\partial u}.
\end{equation}
In the second form of the VP, the differentiability follows from the observation that
$${\cal F}[u\circ\phi_\epsilon]=\int_DL(u,(\nabla\phi_\epsilon)\circ\phi_{-\epsilon}\nabla u)\,\frac{dy}{J(\phi_\epsilon)}\,,$$
where $\phi_{-\epsilon}$ is the diffeomorphism inverse to $\phi_\epsilon$, and $J$ denotes the Jacobian determinant. Here we have used the change of variables $y\mapsto\phi_\epsilon(y)$. Because $\phi$ is a smooth function of both $y$ and $\epsilon$, and thanks to the expansions
\begin{equation}
\label{eq:exp}
J(\phi_\epsilon)=1+\epsilon{\rm div}\,\xi+O(\epsilon^2),\qquad(\nabla\phi_\epsilon)\circ\phi_{-\epsilon}=I_d+\epsilon\nabla\xi+O(\epsilon^2),
\end{equation}
this yields the {\em second form} of the Euler--Lagrange equations:
\begin{equation}
\label{eq:sfEL}
{\rm Div}\left(\nabla u\otimes\frac{\partial L}{\partial p}\right)=\nabla L.
\end{equation}
An important difference between (\ref{eq:ffEL}) and (\ref{eq:sfEL}) is that the former is a scalar PDE, while the latter is a system. Also, the second is conservative, in the sense that it reads ${\rm Div}\,T=0$ where the tensor $T$ is given by
$$T(x)=\left(p\otimes\frac{\partial L}{\partial p}\right)(u,\nabla u)-L(u,\nabla u)\,I_d.$$

When $u$ is of class $C^1$, the second form follows from the first one, by multiplying (\ref{eq:ffEL}) by $\nabla u$ and applying the chain rule. But if $\nabla u$ displays discontinuities instead, the jump relations (Rankine--Hugoniot relations) given by (\ref{eq:ffEL}) and (\ref{eq:sfEL}) are independent from each other. They are respectively
\begin{equation}
\label{eq:ffRH}
\left[\frac{\partial L}{\partial p}\cdot\vec\nu\right]=0
\end{equation}
and
\begin{equation}
\label{eq:sfRH}
\left[T\vec\nu\right]=0,
\end{equation}
where $\vec\nu$ is the unit normal to the discontinuity locus.

Notice that if we were minimizing $\cal F$, then both the first and the second forms of the VP should apply. Thus the likely discontinuities of local minimizers must satisfy both of (\ref{eq:ffRH}) and (\ref{eq:sfRH}). This is often interpreted by saying that the graph of $L(u,\cdot)$ admits a doubly tangent hyperplane. However, models of Mathematical Physics\footnote{We limit ourselves to the study of reversible or quasi-reversible models.} are usually more involved because of the presence of both time and spacial variables. Motions are critical points in some sense, but they are by no means minimizers. There is no reason in general why both forms of the VP should be in force. This will be particularly true in the sequel, when the unknown is a differential form of degree two or more. 

We thus have to choose which form of the VP is meaningful from a physical point of view. The right choice can be guessed if we remember that  Physics is consistent with conservation of energy and momentum: we should privilegiate the second form of the Euler--Lagrange equations, because it gives us the same number of conservations laws as the space-time dimension. Then the analogue of $T$ will, possibly up to the multiplication by a fundamental constant matrix, be the {\em energy-momentum tensor}. The fact that this tensor is always symmetric will follow from a specific property of the Lagrangian density, its invariance with respect to the action of the natural underlying Lie group, its group of symmetries. We remark in passing that in this description, the basic conservation laws (energy, momentum) are those associated with the translation invariance of the Lagrangian, according to Noether's theorem. The symmetry of the tensor $T$ will derive from the `rotational' invariance , where the `rotations' are specific to each model.

Let us give immediately an example within the context of ordinary CV. Here we may assume that $y\in D$ is just a space variable and our VP is actually a minimization. Then isotropy means that $\cal F$ is invariant under $u\mapsto u\circ R$ for every orthogonal transformation $R\in O_d(\R)$. In other words $L(u,Rp)\equiv L(u,p)$ and we can write $L(u,p)=\ell(u,|p|)$, so that our tensor
$$T=\frac{\partial\ell}{\partial r}(u,\nabla u)\,\frac{\nabla u\otimes \nabla u}{|\nabla u|}-L(u,\nabla u)\,I_d$$
is symmetric. Conversely, if the tensor
$$p\otimes\frac{\partial L}{\partial p}-L\,I_d$$ is always symmetric, then $\frac{\partial L}{\partial p}$ is parallel to $p$, and thus $L$ must be isotropic. An example of the situation above is that of the so-called {\em non-parametric} minimal surfaces, where $L(p)=\sqrt{1+|p|^2\,}\,$.

\paragraph{Remark.}
When $L=L(p)$ only, and if $D$ is simply connected, then the unknown $\nabla u$ can be viewed equivalently as a closed differential form $\alpha=du$ of degree one. This point of view inserts the classical Calculus of Variations into the more general formalism described in Section \ref{s:DF}. We notice that $\nabla (u\circ \phi)$ becomes the pullback $\phi^*\alpha$, because of the formul\ae
$$d(\phi^*u)=\phi^*(du),\qquad\phi^*u=u\circ\phi.$$

\paragraph{Plan of the paper.} Section \ref{s:DF} extends the observations above to the situation where the argument of the Lagrangian density $L$ is a closed form of arbitrary degree. The relevant definition of a critical point is that of the {\em second Variational Principle}. The resulting Divergence-free tensor $T$ is described in Theorem \ref{th:main}. We point out that we do not need to introduce a Lagrange multiplier, which is deprived of physical meaning, and whose existence cannot be rigorously proved.
A symmetry property for $T$ is shown to be equivalent to the invariance of $L$ under the action of an orthogonal group (Theorem \ref{th:symvsinv}). Section \ref{s:Phys} displays applications to classical and relativistic gas dynamics, where the forms have degree $n $ in space-time dimension $1+n$, and to the Maxwell system, where the forms have degree $2$.

\paragraph{Warning.} The first form of the VP is irrelevant for our purpose, therefore our analysis is completely different from that done by Bandyopadhyay \& al.  \cite{BDS}.

\section{Calculus of Variations for differential forms}\label{s:DF}

Let $p\in[\![1,d-1]\!]$ be given. We are interested in functionals $\cal F$ whose argument is a closed $p$-form $\alpha$~:
$$d\alpha=0.$$
Mind that the open domain $D\subset\R^d$ may have an arbitrary topology. In particular $\alpha$ does not need to be exact.

An element $\alpha$ of $\Lambda^p(\R^d)$ can be decomposed over the canonical basis
\begin{equation}
\label{eq:ambigu}
\alpha=\sum_{|J|=p}A_Jdx_J,
\end{equation}
where $J=(j_1,\ldots,j_p)$ is ordered by $1\le j_1<\cdots<j_p\le d$, and $dx_J=dx_{j_1}\wedge\cdots\wedge dx_{j_p}$. The sum consists of at most $\binom dp$ terms.
When $\alpha$ is a $p$-form, then the coordinates are functions $y\mapsto A(y)$. The functional we are interested in is 
$${\cal F}[\alpha]=\int_DL(A(y))\,dy,$$
where $L$ is a given smooth function.

\begin{defin}
\label{d:sfEL}
The second Variational Principle for $\cal F$ is  
$$\left.\frac{d}{d\epsilon}\right|_{\epsilon=0}{\cal F}[\phi_\epsilon^*\alpha]=0$$
for every flow associated with a test vector field $\xi\in{\cal D}(D)^d$.
\end{defin}

As in Section \ref{s:clas}, the derivative makes sense under mild assumptions because the change of variables $z=\phi_\epsilon(y)$ yields
$${\cal F}[\phi_\epsilon^*\alpha]=\int_DL((\phi_\epsilon^*\alpha)\circ\phi_{-\epsilon})\,\frac{dz}{J(\phi_\epsilon)}\,,$$
where the argument must be read component-wise.
Our VP thus writes
$$\int_D\left(\frac{\partial L}{\partial A}\cdot\left.\frac{d}{d\epsilon}\right|_{\epsilon=0}((\phi_\epsilon^*\alpha)\circ\phi_{-\epsilon})-L(A){\rm div}\,\xi\right)\,dy=0,\qquad\forall\xi\in{\cal D}(D)^d.$$

\bigskip

Recall that for an arbitrary $\phi\in{\rm Diff}(D)$, the pullback of $\alpha$ is defined by
$$\phi^*\alpha(z)=\sum_IA_I(\phi(z))d\phi_{i_1}\wedge\cdots\wedge d\phi_{i_p}$$
and
\begin{eqnarray*}
(\phi^*\alpha)(\phi^{-1}(y)) & = & \sum_IA_I(y)(d\phi_{i_1}\circ\phi^{-1})\wedge\cdots\wedge (d\phi_{i_p}\circ\phi^{-1}) \\
& = & \sum_{I,J}A_I(y)\left(\det\frac{\partial\phi_I}{\partial y_J}\circ\phi^{-1}\right)\,dy_J.
\end{eqnarray*}
If $\phi_\epsilon$ is a flow at time $\epsilon$, then 
$$(\phi_\epsilon^*\alpha)(\phi_{-\epsilon}(y))=\sum_JA_{J\epsilon}(y)dy_J,\qquad
A_{J\epsilon}(y):=\sum_{I}A_I(y)\left(\det\frac{\partial\phi_{\epsilon I}}{\partial y_J}(\phi_{-\epsilon}(y))\right)\,.$$
Using again (\ref{eq:exp}), there comes
$$A_{J\epsilon}(y)=A_J(y)+\epsilon A_J(y)\sum_{j\in J}\partial_j\xi_j+\epsilon\sum_{d(I,J)=1}(-1)^qA_I(y)\partial_j\xi_i+O(\epsilon^2).$$
The last sum above runs over the $p$-subsets $I$ which differ from $J$ by only one element: $I=K\cup\{i\}$ and $J=K\cup\{j\}$ with $K=I\cap J$. The sign $(-1)^q$ involves the number $q$ of indices $k\in K$ that lie between $i$ and $j$. The whole thing can be recast as
$$A_{J\epsilon}(y)=A_J(y)+\epsilon\sum_{K\subset J}\sum_{i\not\in K}(-1)^qA_I(y)\partial_j\xi_i+O(\epsilon^2)$$
where now $K$ runs over $(p-1)$-subsets of $J$, with again $I:=K\cup\{i\}$ and $\{j\}:= J\setminus K$.
The second VP thus writes
$$\int_D\left(\sum_K\sum_{i,j\not\in K}(-1)^qA_I(y)\frac{\partial L}{\partial A_J}\,\partial_j\xi_i-L(A){\rm div}\,\xi\right)\,dy=0,\qquad\forall\xi\in{\cal D}(D)^d.$$
This yields the second form of the Euler--Lagrange equations:
\begin{thm}
\label{th:main}
A $p$-form $\alpha$ is a critical point of $\cal F$ in the sense of the second Variational Principle, if and only if the tensor 
\begin{equation}
\label{eq:mainT}
T_{ij}  :=  L(A(y))\,\delta_i^j-\sum_{i,j\not\ni K}(-1)^qA_I(y)\frac{\partial L}{\partial A_J}(A(y))
\end{equation}
is Divergence-free:
$$\sum_i\partial_j T_{ij} = 0,\qquad \forall i\in[\![1,d]\!].$$
Hereabove the sum runs over $(p-1)$-subsets $K\subset[\![1,d]\!]$, while $I:=K\cup\{i\}$ and $J:=K\cup\{j\}$ are written as increasing $p$-tuples. The integer $q=q(i,j;K)$ is the number of indices $k\in K$ that lie between $i$ and $j$.
\end{thm}

\paragraph{Changing the order in tuples.}
 In order to keep the calculation at the lowest level of complexity, it can be useful to re-order here and then our $p$-tuples. Let $J$ be an increasing $p$-tuple. If $\sigma\in{\frak S}_p$, then $H=\sigma\cdot J$ represents the same subset of $[\![1,d]\!]$, but $dy_H$ differs from $dy_J$ by the factor $\varepsilon(\sigma)$, the signature. Defining $A_H=\varepsilon(\sigma)A_J$, we ensure that $A_Hdy_H=A_Jdy_J$. For instance, if $p=2$ then $A_{21}=-A_{12}$. Thus we may write equivalently
$$\alpha=\sum A_Hdy_H,$$
where $H$ runs over a collection of $p$-tuples representing the $p$-subsets  of $[\![1,d]\!]$. We emphasize that such a sum is free of repetition. For instance, if $p=2$, then the sum takes in account either $(12)$ or $(21)$, but not both pairs.
Under this convention, the arguments of $L$ may be indexed by any choice of representents of $p$-subsets. We warn however that changing $J$ into $H=\sigma\cdot J$ leads to the identity
$$\frac{\partial L}{\partial A_H}=\varepsilon(\sigma)\,\frac{\partial L}{\partial A_J}\,.$$

An equivalent form of (\ref{eq:mainT}), in which the signs cancel, is obtained by relaxing the condition that the tuples be ordered increasingly:
\begin{equation}
\label{eq:TijK}
T_{ij}=L(A)\,\delta_i^j-\sum_{i,j\not\ni K}A_{iK}\frac{\partial L}{\partial A_{jK}},
\end{equation}
where the sum runs over representents $K$ of $(p-1)$-subsets, and $iK$ denotes the $p$-tuple $(i,k_1,\ldots,k_{p-1})$ obtained by concatenation. 

\paragraph{Symmetry {\em vs} invariance.} 
Let $S\in{\bf Sym}_d$ be a non-degenerate matrix and ${\bf O}(S)$ be its orthogonal group. Although applications to the real word involve either a positive definite matrix, or that of a Minkowski metric, the signature of $S$ is arbitrary in the following analysis.

Matrices $M\in{\bf M}_d(\R)$ act in the natural way over the exterior algebra by pullback. For instance $M^*dy_i=\sum_jm_{ij}dy_j$ and more generally
$$M^*dy_I=\sum_JM\binom IJ\,dy_J,$$
where we use the standard notation for minors. 

Suppose that the Lagrangian density is invariant under the pullback action of the connected component $G$ of ${\bf O}(S)$. Because the exponential map $ \mathfrak{g}\rightarrow G$ defined over the Lie algebra of ${\bf O}(S)$ is onto, this is equivalent to saying that 
\begin{equation}
\label{eq:LNst}
\left.\frac{d}{dt}\right|_{t=0} L\left((e^{tN})^*\alpha\right)=0,\qquad\forall\alpha\in\Lambda^p(\R^d),\,\forall N\in \mathfrak{g}.
\end{equation}
With the notations above, we have
$$\left.\frac{d}{dt}\right|_{t=0} (e^{tN})^*\alpha=\sum_I\sum_{K\subset I}\sum_{j\not\in K}(-1)^qn_{ij}A_Jdy_I=\sum_K\sum_{i,j\not\in K}n_{ij}A_{jK}dy_{iK},$$
so that (\ref{eq:LNst}) rewrites
$$\sum_K\sum_{i,j\not\in K}n_{ij}A_{jK}\frac{\partial L}{\partial A_{iK}}=0,\qquad\forall\alpha\in\Lambda^p(\R^d),\,\forall N\in \mathfrak{g}.$$
In terms of our tensor, we recast this as
$$\Tr(N(L\,I_d-T^T))=0,\qquad\forall\alpha\in\Lambda^p(\R^d),\,\forall N\in \mathfrak{g}.$$
Recalling that
$${\bf O}(S)=\{M\in{\bf M}_d(\R)\,|\,M^TSM=S\},\qquad\mathfrak{g}=\{N\in{\bf M}_d(\R)\,|\,N^TS+SN=0_d\},$$
the condition above is equivalent to
$$\Tr(S^{-1}A(L\,I_d-T^T))=0$$
for every $\alpha$ and every skew-symmetric matrix $A$. This amounts to saying that the corrected tensor $T^TS^{-1}$ is symmetric. We have thus proved 
\begin{thm}
\label{th:symvsinv}
Let $S\in{\bf Sym}_d$ be non-degenerate. Then the tensor $S^{-1}T$ (which is Divergence-free for critical points of the second VP) is symmetric if, and only if the Lagrangian density $L$ is invariant under the pullback action of the neutral component $G$ of the orthogonal group ${\bf O}(S)$.
\end{thm}

Because of the statement above, it is natural to call $S^{-1}T$ the {\em energy-momentum} tensor when the VP describes a physical model. A special case of symmetry-{\em vs}-invariance has been observed for a long time in non-linear models of electro-magnetism --~see for instance \cite{Ser_rel},~-- where $d=4$ and $S$ defines the Minkowski matrix~; see Paragraph \ref{ss:Max}. Another interesting application occurs in relativistic gas dynamics, see Paragraph \ref{ss:RGD}.

\paragraph{Adding a scalar field}
A more general context occurs when $\cal F$ depends upon several forms of various degrees. This happens for instance if we mix two physical aspects, say matter and electro-magnetic field. It happens also when the motion is only quasi-reversible, flows being in general irreversible once shock waves develop. Then a scalar field must be added, namely the entropy $s$. The Lagrangian is now
$${\cal F}[\alpha,s]=\int_DL(A(y),s(y))\,dy,$$
where $s$ must be understood as a differential form of degree zero, though not a closed one~! Definition \ref{d:sfEL} is modified in the obvious way: the second VP writes 
$$\left.\frac{d}{d\epsilon}\right|_{\epsilon=0}{\cal F}[\phi_\epsilon^*\alpha,\phi_\epsilon^*s]=0$$
for every flow associated with a test vector field $\xi\in{\cal D}(D)^d$. We point out that $\phi_\epsilon^*s\equiv s\circ\phi_\epsilon$. An important remark is that the second form of the Euler--Lagrange equations is exactly the same as in Theorem \ref{th:main}, because the change of variable $y\mapsto\phi_\epsilon(y)$ gives back $s$ in the second argument in the integral~; thus the differentiation with respect to $\epsilon$ does not involve the derivative $\partial L/\partial s$ at all.

\paragraph{Comparing with the first VP.} In the first form of the variational principle instead, we write
\begin{equation}
\label{eq:ddeps}
\left.\frac{d}{d\epsilon}\right|_{\epsilon=0}{\cal F}[\alpha+\epsilon d\beta]=0,
\end{equation}
for every compactly supported differential form $\beta$ of degree $p-1$. The corresponding {\em first form} of the Euler--Lagrange equation is
\begin{equation}
\label{eq:dstar}
d^*\left(\frac{\partial L}{\partial \alpha}\right)=0,\qquad\frac{\partial L}{\partial \alpha}:=\sum_J\frac{\partial L}{\partial A_J}\,dy_J.
\end{equation}
This equation takes place in the space of $(p-1)$-forms.

When $\alpha$ is $C^1$-smooth, then (\ref{eq:dstar}) implies ${\rm Div}\,T=0$, because we may write
$$\phi_\epsilon^*\alpha=\alpha+\epsilon{\cal L}_\xi\alpha+O(\epsilon^2),$$
where ${\cal L}_\xi$ stands for the Lie derivative. Considering a $p-1$-form $\gamma$ (possibly multivalued if $D$ has a non-trivial topology) such that $d\gamma=\alpha$, the second VP rewrites therefore as the identity (\ref{eq:ddeps}) for $\beta={\cal L}_\xi\gamma$, which is univalent.

If $\binom{d}{p-1}$ is larger than $d$ then the converse does not hold: the second Euler-Lagrange equations ${\rm Div}\,T=0$ do not imply the first ones (\ref{eq:dstar}), the latter being stronger than the former. We shall illustrate this discrepancy with gas dynamics, where $p=d-1$~: the second VP yields as expected full gas dynamics, while flows obeying to the first VP are actually irrotational.

\section{Applications to time-dependent Physics}\label{s:Phys}

When the Variational Principle deals with time-dependent Physics, one has $d=1+n$ where $n$ is the space dimension. The indices range from $0$ to $n$ and the write $y=(t,x)$ where $t$ is a time variable. The line of index $j=0$ of ${\rm Div}_{t,x}T=0$ stands for the conservation of energy, while the other lines express the conservation of momentum.

\subsection{The case of $n$-forms}

When $p=n=d-1$, a closed $n$-form
$$\alpha=\rho dx_1\wedge\cdots\wedge dx_n+\cdots+(-1)^nq_ndt\wedge dx_1\wedge\cdots\wedge dx_{n-1}$$
corresponds to a scalar conservation law
\begin{equation}
\label{eq:cofm}
\partial_t\rho+{\rm div}_xq=0.
\end{equation}
This is often interpreted as the conservation of mass. In Classical Mechanics $\rho$ is the mass density and $q$ the linear momentum. For the sake of simplicity, we set $m_0:=\rho$ and $m_i=q_i$ otherwise. Denoting $\hat\imath$ the increasing $n$-tuple $(\ldots,i-1,i+1,\ldots)$, we have
$$A_{\hat\imath}=(-1)^im_i.$$

The Lagrangian density being a function $L(m,s)$ with $s$ the entropy scalar field, the sum in (\ref{eq:mainT}) is rather simple. If $i\ne j$, it runs over $(d-2)$-subsets $K$ such that $i,j\not\in K$, and there is only one, namely $K=[\![0,n]\!]\setminus\{i,j\}$. Thus
$$T_{ij}=-(-1)^qA_{\hat j}\frac{\partial L}{\partial A_{\hat\imath}}\,=(-1)^{j-i}A_{\hat j}\frac{\partial L}{\partial A_{\hat\imath}}\,=m_j\frac{\partial L}{\partial m_i}\,,\qquad\forall i\ne j.$$
The diagonal entries are instead given by
$$T_{ii}=L-\sum_{k\ne i}A_{\hat k}\frac{\partial L}{\partial A_{\hat k}}=L(m)+m_i\frac{\partial L}{\partial m_i}-m\cdot\frac{\partial L}{\partial m}.$$
Overall we obtain
\begin{equation}
\label{eq:Tnform}
T=\frac{\partial L}{\partial m}\otimes m+\left(L-m\cdot\frac{\partial L}{\partial m}\right)\,I_d.
\end{equation}

\paragraph{Transport of the entropy.}
If the flow is smooth enough so as the chain rule be valid, then we have
$$0={\rm div}\,T_{i\bullet}=({\rm div}\,m)\frac{\partial L}{\partial m_i}+m\cdot\nabla\frac{\partial L}{\partial m_i}-m\cdot\partial_i\frac{\partial L}{\partial m}+\frac\partial{\partial s}\left(L-m\cdot\frac{\partial L}{\partial m}\right)\,\partial_is,$$
with the first term vanishing identically because of $d\alpha=0$. Multiplying by $m_i$ and then summing over $[\![0,n]\!]$, there remains
$$\frac\partial{\partial s}\left(L-m\cdot\frac{\partial L}{\partial m}\right)\,m\cdot\nabla s=0.$$
Provided that the function
$$\frac\partial{\partial s}\left(L-m\cdot\frac{\partial L}{\partial m}\right)$$
does not vanishes, we infer $m\cdot\nabla s\equiv0$, which is the transport of the entropy along the particle paths. Of course this calculation is not valid any more across discontinuities.

\subsubsection{Gas dynamics}\label{ss:GD}

In classical mechanics, the Lagrangian density is the difference of the kinetic energy and of a function of so-called internal variables:
\begin{equation}
\label{eq:kinmsint}
L(\rho,q,s)=\frac{|q|^2}{2\rho}-g(\rho,s).
\end{equation}
We notice that the kinetic term, being homogeneous of order one, does not contribute to the scalar part of $T$ (Euler identity): we just have
$$L-m\cdot\frac{\partial L}{\partial m}=\rho\partial_\rho g-g.$$ 
With
\begin{equation}
\label{eq:dLdm}
\frac{\partial L}{\partial m}=\begin{pmatrix}
-\frac{|q|^2}{2\rho^2}-\partial_\rho g \\ \frac q\rho 
\end{pmatrix},
\end{equation}
we obtain
$$T=\begin{pmatrix}
-\frac{|q|^2}{2\rho}-g &  -\left(\frac{|q|^2}{2\rho^2}+\partial_\rho g\right)q^T \\ q & \frac{q\otimes q}\rho+(\rho\partial_\rho g-g)I_n
\end{pmatrix}\,.$$
That each row is divergence-free expresses the conservation of energy (first row) and of momentum (the other ones). The pressure is given by
$$p(\rho,s):=\rho\partial_\rho g-g,$$
while $g$ is the internal energy per unit volume.

\paragraph{Symmetry.}
Because the model is invariant under the action of ${\bf O}(n)$, but not of any orthogonal group of $\R^{1+n}$, our tensor $T$ cannot be made symmetric by multiplying by a constant matrix. Only the lower-right block is symmetric. To recover a full symmetry, we must take into account the constraint (\ref{eq:cofm}) and replace the first row of $T$ by the vector field $m$~:
$$T'=\begin{pmatrix}
\rho & q^T \\ q & \frac{q\otimes q}\rho+(\rho\partial_\rho g-g)I_n
\end{pmatrix}\,.$$

Conversely, in order that the tensor obtained by replacing the first row of $T$ by $m$ be symmetric, we need
$$\rho\frac{\partial L}{\partial q}=q,$$
which implies an expression of the form (\ref{eq:kinmsint}). This shows that the symmetry of $T'$ is tightly related to the specific choice of the Lagrangian, in particular to the form of the kinetic energy.

\paragraph{Remark.} Our approach by the second VP  is essentially equivalent to the calculations done in \cite{Ser_M2AN}.

\paragraph{Flows obeying the first VP.}
If instead we perform variations $\alpha\mapsto\alpha+\epsilon d\beta$ at the target --~the first Variational Principle,~-- we obtain the first form (\ref{eq:dstar}) of the Euler--Lagrange equations. With (\ref{eq:dLdm}), this means on the one hand that the velocity field $v:=\frac q\rho$ is irrotational, and on the other hand that the Bernoulli relation holds true:
$${\rm curl}_x\,v=0,\qquad\partial_tv+\nabla_x\left(\frac{|v|^2}2+\partial_\rho g\right)=0.$$
As usual, this can be recast by introducing a velocity potential $\psi$~:
$$v=\nabla_x\psi,\qquad\partial_t\psi+\frac12\,|\nabla_x\psi|^2+\partial_\rho g=0.$$

\subsubsection{Relativistic gas dynamics}\label{ss:RGD}

For the sake of simplicity, we shall focus on special relativity, where $c$ is the light speed and $n=3$. The space $\R^{1+3}$ is equipped with the Minkowski metric $-c^2dt^2+dx_1^2+dx_2^2+dx_3^2$ associated with the matrix $\Lambda={\rm diag}(-c^2,1,1,1)$. The motion of the flow is described by the divergence-free $4$-momentum vector $m(t,x)$ which satisfies the constraint
$$m_0\ge\frac1c\,\sqrt{m_1^2+m_2^2+m_3^2\,}\,.$$
The particle number density is defined by
$$\rho:=\sqrt{-m^T\Lambda m\,}=\sqrt{c^2m_0^2-m_1^2-m_2^2-m_3^2\,}\,.$$
The $4$-velocity vector $u=\frac m\rho$ satisfies $u^T\Lambda u\equiv-1$. Notice that we depart slightly from the notations of \cite{MTW}.

\bigskip

The Lagrangian density is a function of $\rho$ and of the entropy, $L=L(\rho,s)$. According to  (\ref{eq:Tnform}), our Divergence-free tensor writes
$$T=\rho\frac{\partial L}{\partial\rho}\,\begin{pmatrix} 
c^2u_0^2 & c^2u_0v^T \\ -u_0v & -v\otimes v
\end{pmatrix}+\left(L-\rho\frac{\partial L}{\partial\rho}\right)\,I_4.$$
As anticipated by Theorem \ref{th:symvsinv}, the symmetry is recovered by multiplying by $-\Lambda^{-1}$, obtaining the (still Divergence-free) energy-momentum tensor
$$T':=-\Lambda^{-1}T=\rho\frac{\partial L}{\partial\rho}\,u\otimes u+\left(\rho\frac{\partial L}{\partial\rho}-L\right)\,\Lambda^{-1}.$$
Introducing the energy density and the pressure by
$$e:=\frac1{c^2}\,L,\qquad p:=\rho\frac{\partial L}{\partial\rho}-L,$$
 this becomes the well-known formulation
$$T'=(ec^2+p)u\otimes u+p\Lambda^{-1}.$$

When $L=f(s)\rho^\kappa$ is homogeneous of degree $\kappa>1$ in the particle density, then $p$ is a linear function of $e$~:
$$p=(\kappa-1)e c^2.$$
The hypothesis that the sound waves travel at sub-luminous velocity is equivalent to $\kappa<2$. Ultrarelativistic jets are believed to correspond to $\kappa=\frac43\,$, see \cite{Ani}.

\bigskip

Let us mention in passing the {\bf limit case} where $L\equiv\rho^2$. Then the seamingly overdetermined quasilinear system
\begin{eqnarray*}
{\rm div}_{t,x}m & = & 0, \\
{\rm Div}_{t,x}\left(2m\otimes m+\rho^2\Lambda^{-1}\right) & = & 0
\end{eqnarray*}
is actually consistent and its waves travel at light speed in every direction. This is confirmed by the study of the Rankine--Hugoniot condition accross a hypersurface with $4$-normal vector $\vec\nu$. The jumps satisfy
$$[m\cdot\vec\nu]=0,\qquad2[m(m\cdot\vec\nu)]+[\rho^2]\Lambda^{-1}\vec\nu=0.$$
Multiplying the second identity by $\vec\nu$ and using the first one to eliminate the jump of $(m\cdot\vec\nu)^2$, there remains
$$[\rho^2]\vec\nu^T\Lambda^{-1}\vec\nu=0.$$
Since $[\rho]\ne0$ accross a genuine discontinuity, we obtain that $\vec\nu^T\Lambda^{-1}\vec\nu=0$. This means that the discontinuity locus propagates at light speed in its normal direction. Amazingly, this is the same law as for the linear wave equation.

\subsection{Calculus for $2$-forms: Maxwell's equations}\label{ss:Max}

When $p=2$, (\ref{eq:TijK}) writes
\begin{equation}
\label{eq:Tkij}
T_{ij}=L(A)\,\delta_i^j-\sum_{k\neq i,j}A_{ik}(y)\frac{\partial L}{\partial A_{jk}}.
\end{equation}
The electromagnetic field in vacuum and Minkowski metric $-dt^2+dx_1^2+dx_2^2+dx_3^2$ turns out to be a closed $2$-form
$$\alpha=(\vec E\cdot dx)\wedge dt+B_1dx_2\wedge dx_3+B_2dx_3\wedge dx_1+B_3dx_1\wedge dx_2.$$
The closedness $d\alpha=0$ translates as the Gau\ss--Faraday law
$$\partial_t\vec B+{\rm curl}_x\vec E=0,\qquad {\rm div}_x\vec B=0.$$
We thus identify
$$A_{j0}=E_j,\quad A_{ij}=\epsilon(ijk)B_k,$$
where $\epsilon$ denotes the signature in $\frak S_3$.

\bigskip

Given a Lagrangian density $L(\vec E,\vec B,s)$, we define the auxiliary vector fields
$$\vec D=\frac{\partial L}{\partial\vec E}\,,\quad \vec H=-\frac{\partial L}{\partial\vec B}\,.$$
From the formula (\ref{eq:Tkij}) the Divergence-free tensor of the second VP reads
$$T=\begin{pmatrix} L-\vec E\cdot\vec D & \vec H\times\vec E \\
\vec D\times\vec B & (L+\vec B\cdot\vec H)\,I_3-\vec E\otimes\vec D-\vec H\otimes\vec B
\end{pmatrix}.$$
The (divergence-free) first row expresses the conservation of energy (Poynting's Theorem)
$$\partial_tW+{\rm div}_x( \vec E\times\vec H)=0,\qquad W:=\vec E\cdot\vec D -L.$$

\paragraph{Symmetry.}
We recall that \cite{Ser_rel} the (still Divergence-free) tensor
$$\tilde T:={\rm diag}(-1,1,1,1)T$$
is symmetric if, and only if $L$ is a function of
$$\frac12(|\vec E|^2-|\vec B|^2),\qquad\vec E\cdot\vec B.$$
Once again, this  illustrates Theorem \ref{th:symvsinv}, since this means that $L$ is Lorentz invariant -- equivalently, that ${\cal F}[\alpha]$ does not depend upon the choice of linear coordinates but only upon the class of $\alpha$ modulo the linear transformations that preserve the Minkowski metric. 

Notice that $\vec E\cdot\vec B$ is nothing but the Pfaffian of $\alpha$. When $L$ reduces to $\frac12(|\vec E|^2-|\vec B|^2),$ the first form of the Euler--Lagrange equation is the Maxwell--Amp\`ere equation~; with the constraint $d\alpha=0$, this gives us back the usual linear Maxwell's equations.




\begin{thebibliography}{00}





\bibitem{Ani} A. M. Anile. {\em Relativistic fluids and magneto-fluids}. Cambridge Monographs on mathematical Physics. Cambridge Univ. Press (1989).



\bibitem{BDS} S. Bandyopadhyay, B. Dacorogna, S. Sil. Calculus of variations with differential forms. {\em J. Eur. Math. Soc.}, {\bf17} (2015), pp 1009--1039.


\bibitem{MTW} C. W. Misner, K. S. Thorne, J. A. Wheeler. {\em Gravitation} (1973) W. H. Freeman \& Co (San Francisco). 


\bibitem{Ser_M2AN} D. Serre. Sur le principe variationnel des \'equations de la m\'ecanique des fluides parfaits. {\em M2AN - Mod\'el. Math. et Anal. Num.}, {\bf27} (1993), pp 739--758.


\bibitem{Ser_rel} D. Serre. Non-linear electromagnetism and special relativity. {\em Discrete Cont. Dynam. Syst.}, {\bf23} (2009), pp 435--454.

\bibitem{Ser_DPT} D. Serre. Divergence-free positive symmetric tensors  and fluid dynamics. {\em Annales de l'Institut Henri Poincar\'e (analyse non lin\'eaire)}, {\bf35} (2018), pp 1209--1234.

\bibitem{Ser_JMPA} D. Serre. Compensated integrability. Applications to the Vlasov--Poisson equation and other models of mathematical physics. {\em J. Math. Pures \& Appl.}, {\bf127} (2019), pp 67--88.





 

\end{thebibliography}
\end{document}